%% file: main-v2.tex
\documentclass[11pt]{article} 

\usepackage[utf8]{inputenc} 

\usepackage{cancel}

\usepackage{geometry} 
\usepackage{graphicx} 
\usepackage[english]{babel}

\usepackage{booktabs} 
\usepackage{array} 
\usepackage{paralist} 
\usepackage{verbatim} 
\usepackage{subfig} 
\usepackage{amssymb}
\usepackage{bbm}

\usepackage{amsmath}
\usepackage{amsthm}
\usepackage{makecell}

\usepackage{epigraph}

\usepackage{hyperref}
\hypersetup{
    colorlinks=true,
    citecolor=blue
}

\usepackage{pdflscape}
\usepackage{natbib}
\usepackage{algorithm2e}

\usepackage{wrapfig}
\usepackage{subfig}
\usepackage[capitalize,noabbrev]{cleveref}
\usepackage{authblk}
\usepackage{changepage}

\usepackage{tikz-cd}
\usepackage{xcolor}

\input{macros}

\newtheorem{theorem}{Theorem}[section]
\newtheorem{definition}[theorem]{Definition}

\newtheorem{fact}{Fact}
\newtheorem{proposition}[theorem]{Proposition}

\newtheorem{lemma}[theorem]{Lemma}
\theoremstyle{definition}

\Crefname{fact}{Fact}{Facts}
\Crefname{assumption}{Assumption}{Assumptions}

\title{Testing by betting while borrowing and bargaining
}

\author[1]{Hongjian Wang}
\author[2]{Wouter M.\ Koolen}
\author[3]{Aaditya Ramdas}
\affil[1, 3]{Carnegie Mellon University}
\affil[2]{University of Twente and CWI}

\date{\today}

\begin{document}

\maketitle

\begin{abstract}
   {
  Testing by betting has been a cornerstone of the game-theoretic statistics literature. One bets against the null hypothesis, and the accumulated wealth $W_t$ quantifies the evidence against the null hypothesis after $t$ rounds, and the null can be rejected at level $\alpha$ whenever  $W_t \geq 1/\alpha$.
  A key assumption permeating the literature is that one cannot bet more money than they currently have (the wealth must stay nonnegative). 
  In this work, we examine the consequences of allowing the bettor to borrow money in each round (for example after going bankrupt). Specifically, we ask how the threshold of $1/\alpha$ must be accordingly adjusted to retain the desired level $\alpha$.
  Our findings are twofold. First, if the new rejection rule is $W_t \geq g(\alpha,L_t)$ where $L_t$ is the total liability at time $t$, then we show that $g(\alpha,0)>1/\alpha$ if $g(\alpha,L_t)<\infty$ for any $L_t > 0$; in words, we must pay for the possibility of borrowing, even if in fact we do not borrow. Second, and in contrast to the first, if one employs a \emph{path dependent threshold} $h(\alpha,W_0,L_1,\dots,W_{t-1},L_t)$, that is a function of past \emph{leverage ratios}, then there is in fact no extra price to pay for the possibility of borrowing.}
\end{abstract}

\section{Introduction}
Aggressive gamblers put borrowed money on the table to increase their potential win. Risk-seeking traders enter leveraged positions for exposure to higher expected returns. In this paper, we study the analogous situations in \emph{testing by betting}, a core foundational principle in game-theoretic statistics \citep{shafer2021testing,ramdas2023game}.

In testing by betting, the statistician reports the outcome of a bet against the null hypothesis, rather than a p-value, as the measure of statistical evidence. Among its numerous benefits, adopting the testing by betting framework naturally allows for a sequential (``opportunistic'') testing procedure as outlined by \citet[Protocol 2, Section 4.1]{shafer2021testing}: the statistician places a bet, observes a variable, receives the payoff from the bet, and decides to repeat or not. In the meantime, the betting statistician's \emph{wealth process} $\{W_t\}_{t \ge 0}$ measures the amount of evidence accrued against the tested null hypothesis, and, importantly, evolves as a martingale under the null distribution. 
As an application of Ville's theorem, the statistician may reject the null when $W_t \ge 1/\alpha$, where $\alpha \in (0,1)$ is a prescribed type 1 error level. We will soon formalize this standard framework in our terms in \cref{sec:noborrow}.

Bettors in real life may borrow and leverage to bet. However, the game-theoretic statistics literature has barred the analogue in testing by betting, insisting that a fundamental bedrock of the testing-by-betting framework is that of \emph{not betting more than one has}, so that the wealth process $\{ W_t \}$ is nonnegative. In this paper, we try to extend the classical testing-by-betting framework to allow for such situations. In particular, we study the types of rejection rules that still allow testing hypotheses at the predefined level $\alpha$, generalizing the use of Ville's theorem in the standard setup summarized above. 


\section{Adding borrowing to testing by betting}
\subsection{Background: Classical betting game without borrowing}\label{sec:noborrow}

Let us start with the standard sequential testing-by-betting setup without borrowing. Two players, Casino\footnote{This player is referred to as \emph{Forecaster} in Shafer's original writings, but we shall later see that it is conceptually much more convenient to use Casino as the analogy in our work.} and Statistician are involved in a game, where Casino offers
bets to Statistician to test the verity of a null hypothesis. Before each {random} outcome is revealed, Statistician places their wealth on the outcome, at odds specified by the null hypothesis. The total wealth of Statistician after one or more rounds of betting can be interpreted as the amount of \emph{evidence} gathered against the null. This is often referred to as ``testing by betting'', and the corresponding sequential testing paradigm usually involves a nonnegative martingale or, more generally, a nonnegative supermartingale (``NSM'' for short) denoting the sequential wealth evolution. We rephrase and slightly generalize the ``testing protocols'' of \citet[Section 4.1]{shafer2021testing} via the following definition.

\begin{definition}[NSM Betting]\label[definition]{def:nsm-bet}  \normalfont
    Let $(\Omega, \cA, \Pr)$ be a probability space with filtration $\{\cF_t\}_{t \ge 0}$. Let $\{B_t \}_{t \ge 1}$ be a nonnegative process adapted to $\{ \cF_t \}_{t \ge 1}$. Statistician's wealth process $\{ W_t \}_{t \ge 0}$ is defined as $W_0=1$ and
\begin{equation}
     W_t = W_{t-1} \cdot  B_t.
\end{equation}
We consider two nested null hypotheses:
\begin{itemize}
    \item If $\Exp(B_t | \cF_{t-1} ) \le 1$ for all $t \ge 1$,  we denote the probability by $\Ps$ and the expected value by $\Exps$. This is called the \emph{supermartingale null}.
    \item Further, if $\Exp(B_t | \cF_{t-1} ) = 1$ for all $t \ge 1$, we denote the probability by $\Ph$ and the expected value by $\Exph$. This is called the \emph{martingale null}.
\end{itemize}
\end{definition}

As a simple economic interpretation of \cref{def:nsm-bet},  $\{B_t\}$ is the random per-unit payoff at time $t$ of Statistician's portfolio allocation at time $t-1$. 
Casino may offer various bets to Statistician, and Statistician may allocate their wealth $W_{t-1}$ arbitrarily between the bets offered (and cash). On aggregate, 1 unit of wealth invested in this portfolio at time $t-1$ becomes $B_t$ units at time $t$; and Statistician's wealth consequently becomes $W_{t-1} \cdot B_t$ at time $t$. The null $\Ph$ then states that the entire portfolio is always fairly priced, under which $\{ W_t \}$ is a nonnegative martingale; The null $\Ps$, on the other hand, states that the portfolio is at least fairly priced (with the possibility of being overpriced, to the disadvantage of Statistician), under which $\{ W_t \}$ is a nonnegative supermartingale.

An illustrative example of \cref{def:nsm-bet} is the following \emph{fair-coin game} \citep[Chapter 3.1]{shafer2005probability}: $X_t = \pm 1 $ are i.i.d.\ random variables, and $\lambda_t \in [-1, 1]$ are predictable w.r.t.\ the natural filtration of $\{ X_t \}$. Let $B_t = 1+\lambda_t X_t$. Then, \cref{def:nsm-bet} corresponds to the game where Casino offers \emph{double-or-nothing} bets on the outcome of coin tosses $\{ X_t \}$ (heads and tails are 1 and -1 respectively), where Statistician bets $|\lambda_t|$ portion of current wealth $W_{t-1}$ on heads (if $\lambda_t > 0$) or tails (if $\lambda_t < 0$) of the outcome $X_t$. Both the martingale null $\Ph$ and the supermartingale null $\Ps$ state that the coin is fair: $\Ph(X_1 = 1) = \Ps(X_1 = 1) = 1/2$. 

Crucially, the process $\{ W_t \}$ defined in \cref{def:nsm-bet}, as a nonnegative supermartingale under $\Ps$,
satisfies \emph{Ville's inequality} \citep{ville1939etude}. For any $\alpha \in (0,1)$,
\begin{equation}\label{eqn:time-unif-nb} \tag{V.I.}
    \Ps \left( \sup_t W_t \ge W_0/\alpha \right) = \Ps \left( \sup_t W_t \ge 1/\alpha \right)  \le \alpha.
\end{equation}
justifying the notion of ``wealth as evidence'': if $\Ps$ were true, it would be unlikely for the wealth $\{ W_t \}$ to ever be large. A large value of $W_t$ at any time means that $\Ps$ is unlikely to be the ground truth. \eqref{eqn:time-unif-nb} also leads to the rejection rule employed by Statistician:
\begin{equation}\label{eqn:reject-rule}
    \text{Reject $\Ps$ (and hence $\Ph$) when $W_t$ first exceeds $1/\alpha$,}
\end{equation}
as \eqref{eqn:time-unif-nb} ensures that this rule has a type 1 error rate at most $\alpha$ for any prescribed level $\alpha$ (say 0.05). This NSM testing framework has been instantiated for a wide range of parametric and nonparametric testing problems, as summarized recently by \citet[Sections 4, 5]{ramdas2023game}.


\subsection{Our contribution: Betting game with borrowing}

In the framework of testing by betting reviewed above, an overarching assumption has been made throughout the literature: Statistician, the bettor, is not allowed to risk more than what they have; the per-round payoff process $\{ B_t \}$ and the wealth process $\{ W_t \}$ must consequently always be nonnegative.
For example, in the prototypical game testing the fairness of a coin, \citet[Section 3.1]{shafer2005probability} allow the bettor to freely choose the volume of bets but \emph{they lose the game immediately if their wealth becomes negative}, essentially banning any bets larger than the current wealth. The same ``losing at bankruptcy" rule applies to the entire book by these authors, as well as in their follow-up book \citep{shafer2019game}. 
This constraint has essentially remained unchanged and unchallenged in the game-theoretic statistics literature~\citep{ramdas2023game}.

In this paper, we challenge this assumption and allow the betting Statistician to risk more than they have. That is, Statistician can \emph{borrow} money before placing each bet {without being forced to quit the game in the possible event of indebtedness}. Formally, the ``testing by betting with borrowing'' protocol we work with throughout this paper involves the following generalization of \cref{def:nsm-bet}.  \begin{definition}[Borrowed NSM Betting]\label[definition]{def:nsm-bet2}
\normalfont
     Let $(\Omega, \cA, \Pr)$ be a probability space with filtration $\{\cF_t\}_{t \ge 0}$. Let $\{B_t \}_{t \ge 1}$ be a nonnegative process adapted to $\{ \cF_t \}_{t \ge 1}$, and $\{ \beta_t \}_{t \ge 1}$ a process predictable with respect to $\{ \cF_t \}_{t \ge 0}$ such that $\beta_t \ge - W_{t-1}$ for each $t$. Statistician's wealth process $\{ W_t \}_{t \ge 0}$ is defined as $W_0=1$ and
\begin{equation}
     W_t = (W_{t-1} + \beta_t) \cdot  B_t.
\end{equation}
The supermartingale null $\Ps$ under which $\Exps(B_t | \cF_{t-1} ) \le 1$ for all $t \ge 1$, and the martingale null $\Ph$ under which $\Exph(B_t | \cF_{t-1} ) = 1$ for all $t \ge 1$ are the same as in \cref{def:nsm-bet}.
\end{definition}
The process $\{ \beta_t \}$ is the borrowed amount each round. In particular, we allow $\beta_t$ to be negative, which could correspond to Statistician ``paying back" earlier debts. The predictability of $\{\beta_t\}$ in contrast to the adaptedness of $\{ B_t \}$ means that in each round, Statistician first borrows, then bets, and finally observes the outcome. Let us also track the total borrowed amount up to time $t$,
\begin{equation}\label{eqn:liabilities}
    L_t = \sum_{i=1}^t \beta_i
\end{equation}
and call $\{ L_t \}$ the \emph{liabilities} process in this game.

\textbf{The central objective} of this paper is to extend the simple time-uniform inequality \eqref{eqn:time-unif-nb} 
as well as its implied rejection rule \eqref{eqn:reject-rule}
 to this borrowed setting. 
Now that $\{W_t\}$ is no longer a supermartingale under $\Ps$, we ask, can we somehow compensate for it by increasing the threshold $1/\alpha$ accordingly? 

We give a pair of interesting answers to the question above. 

\begin{itemize}
    \item First, if we use a rule of the form ``reject when $W_t$ exceeds $g(\alpha, L_t)$'' for some ``threshold function'' $g$, then $g(\alpha, 0) > 1/\alpha$, unless $g(\alpha, x) = \infty$ for almost all positive $x$. This means that if we want a nontrivial threshold for positive liabilities, we must necessarily pay a price on all paths when we do not actually borrow money. That is, \emph{in order to potentially borrow on \underline{some} paths, there is a price to be paid on \underline{all} paths}.
    \item Second, if we use a rule of the form ``reject when $W_t$ exceeds $h(\alpha, W_0, L_1, \dots, W_{t-1}, L_t )$'', where the threshold function $h$ now allows a dynamic dependence on the wealth and borrowing history,
    then it is possible to design rules such that $h(\alpha,w_0,0,\dots, w_{t-1},0)=1/\alpha$, and $h(\alpha, w_0,\ell_1,\dots, w_{t-1},\ell_t ) < \infty$ for any arguments. This means that path-dependent thresholds can recover this apparent price from the previous method, which only accounted for the total amount borrowed.  However, this approach comes with another tradeoff:  going all-in, losing all, and borrowing to start again eliminates all accumulated evidence.
\end{itemize}

The contrast of the two kinds of results outlined above highlights the very bewildering consequences of allowing a sequence of pre-bet borrowing $\{ \beta_t \}$ into the game. We unpack these in \cref{sec:gL} and \cref{sec:current} respectively.


\subsection{Notations}

We work with the probability space as formulated in \cref{def:nsm-bet2}, primarily focusing on the case where the NSM null $\Ps$ is satisfied, and we denote $\Ps$ for probability and $\Exps$ for expected value. Sometimes, however, we shall prove independent lemmas that hold generally, where we use $\Pr$ and $\Exp$ to highlight such independence. The time index $t$ is always discrete, $t = 0,1,\dots$, when we define a sequence, say,
$\{ x_t \}_{t \ge 0}$. Continuous-time stochastic processes, however, \emph{will} occasionally be introduced, but they will be announced case by case.

Following the definition of the liabilities process $\{ L_t \}$ in \eqref{eqn:liabilities}, we also define
the \emph{net wealth} process for the difference between ``gross" wealth and liabilities,
\begin{equation}
    N_t  = W_t-  L_t.
\end{equation}

Recall that in the previous setting without borrowing in \cref{sec:noborrow}, the wealth process is a nonnegative martingale under $\Ph$ {and a nonnegative supermartingale under $\Ps$}. In our current case with borrowing, we have the following easy-to-note fact on the more general relation between these three processes.

\begin{fact}[Doob decomposition of wealth]\label[fact]{prop:doob-wealth}
    The net wealth process $\{ N_t \}$ is a martingale on $\{ \cF_t \}$ under $\Ph$ {and a supermartingale under $\Ps$}. Further, {under $\Ph$,} $W_t = N_t + L_t$ is the Doob decomposition of the process $\{W_t\}$ into the martingale $\{N_t\}$ and the predictable process $\{ L_t \}$. 
\end{fact}

That is, the nonnegative supermartingale property that holds earlier for $\{ W_t \}$ is now broken down into the new gross wealth $\{ W_t \}$ which is nonnegative but not necessarily a supermartingale, and the net wealth process $\{ N_t \}$ which is a supermartingale but not necessarily nonnegative.

\section{Liabilities curve-crossing probability of gross wealth}\label{sec:gL}

First, we explore the following generalization to \eqref{eqn:time-unif-nb}: we seek to establish curve-crossing probability bounds in the form of
\begin{equation}\label{eqn:L-objective}
   \Ps( \exists t, \ W_t \ge g(\alpha, L_t) ) \le \alpha, \quad \forall \alpha \in (0,1)
\end{equation}
where $g$ is some function. In that way, we can use the 
rule ``reject when $W_t$ exceeds $g(\alpha, L_t)$'' to reject the null $\Ps$.

In fact, with the following main result, we can allow arbitrary ``curve shapes''. The only assumption we need is the (weak) monotonicity of the liabilities process, or equivalently, the nonnegativity of borrowings.

\begin{theorem}\label{thm:g-crossing} Assume that $\{L_t\}$ is non-decreasing. That is, $\Ps(\beta_t \ge 0) = 1$ for all $t$.
Then, for any  function $g:[0,\infty) \to [1,\infty)$ such that $\ell \mapsto g(\ell) - \ell$ is increasing,
    \begin{equation}\label{eqn:g-bound}
    \Ps( \exists t, \ W_t \ge g(L_t) ) \le 1 - \frac{g(0) - 1}{g(0)}\exp \left( -\int_0^\infty \frac{1}{g(\ell)} \d \ell \right).
\end{equation}
Additionally, the inequality \eqref{eqn:g-bound} is tight in the sense that, for a fixed $g$,
\begin{equation}\label{eqn:g-bound-tight}
  \sup  \Ps( \exists t, \ W_t \ge g(L_t) ) = 1 - \frac{g(0) - 1}{g(0)}\exp \left( -\int_0^\infty \frac{1}{g(\ell)} \d \ell \right),
\end{equation}
where the supremum is over all possible instances of the borrowed betting game $(\{B_t\}, \{\beta_t\})$.
\end{theorem}
We prove \cref{thm:g-crossing} in \cref{sec:pf-gcross}.
The assumption that $\{ L_t \}$ must not decrease brings up an intriguing discussion in \cref{sec:why-not-move-back}. 
In order to obtain the desideratum \eqref{eqn:L-objective} with \cref{thm:g-crossing}, we can construct a continuum of functions $\{ g_{\alpha}(\cdot) \}_{\alpha \in (0,1)}$ such that
\begin{equation}
     1 - \frac{g_\alpha(0) - 1}{g_\alpha(0)}\exp \left( -\int_0^\infty \frac{1}{g_\alpha(\ell)} \d \ell \right) = \alpha,
\end{equation}
 and set $g(\alpha, \ell) = g_\alpha (\ell)$.

Immediately, however, we see the following ``deficiency'' of \cref{thm:g-crossing} compared to the \eqref{eqn:time-unif-nb} in the borrow-free setting.

\begin{proposition}\label[proposition]{prop:g-deficiency} In \cref{thm:g-crossing},
    if there exists $\alpha \in (0,1)$ such that the RHS of \eqref{eqn:g-bound} equals $\alpha$, then $g(0) \ge 1/\alpha$. Additionally, if $g(0) = 1/\alpha$, then $g(\ell) = \infty$ for all $\ell > 0$. 
    Contrapositively, if $g(\ell)$ is finite for any $\ell > 0$, then $g(0) > 1/\alpha$.
\end{proposition}
\begin{proof}
    Note that
\begin{align}\label{eqn:01expression}
   (1-\alpha) \frac{g(0)}{g(0) - 1} = \exp \left( -\int_0^\infty \frac{1}{g(\ell)} \d \ell \right) \in (0,1]
\end{align}
where the inclusion holds since the integration above is nonnegative. Solving $ (1-\alpha) \frac{g(0)}{g(0) - 1} \le 1$, we obtain $g(0) \ge 1/\alpha$.

Next, if $g(0)  = 1/\alpha$, then the expression in \eqref{eqn:01expression} equals 1, 
so
\begin{equation}
    \int_0^\infty \frac{1}{g(\ell)} \d \ell = 0,
\end{equation}
meaning that $1/g(\ell) = 0$, i.e.\ $g(\ell) = \infty$ almost everywhere. Since $g(\ell) - \ell$ is increasing, we see that $g(\ell) = \infty$ for all $\ell > 0$, concluding the proof. 
\end{proof}

In words, \cref{prop:g-deficiency} states that using the liabilities $L_t$ to adjust the reject threshold via \cref{thm:g-crossing} will always, no matter how ``smart'' the thresholding curve $g()$ is picked, lead to the threshold being above $1/\alpha$ even when Statistician is not borrowing. That is, there is one price $g(L_t) - g(0) > 0$ to pay for \emph{actually} borrowing $L_t$ in the current path, and another general, path-independent price $g(0) - 1/\alpha > 0$ to pay for \emph{potentially} borrowing in any possible path.

Finally, some quadratic examples of $g$ functions are provided in \cref{sec:g-inventory}.

\section{Leverage penalty-crossing probability of gross wealth}\label{sec:current}

Having considered (and seen the deficiency of) the generalization to \eqref{eqn:time-unif-nb} that replaces the $1/\alpha$ threshold with $g(\alpha, L_t)$, we now relax our goal and open up to more flexibility by seeking bounds of form
\begin{equation}\label{eqn:dep-objective}
   \Ps( \exists t, \ W_t \ge h(\alpha, W_0, L_1,\dots, W_{t-1}, L_t) ) \le \alpha, \quad \forall \alpha \in (0,1),
\end{equation}
where the threshold $h(\alpha, \cdots)$ now take the \emph{history} of past liabilities $(L_1, \dots, L_t)$ as well as the past wealths $(W_0, \dots, W_{t-1})$ as arguments. If \eqref{eqn:dep-objective} can be established, we use the 
rule ``reject when $W_t$ exceeds $h(\alpha, W_0, L_1,\dots, W_{t-1}, L_t)$'' to reject the null $\Ps$.

Our result in this section achieves this objective by considering a different notion of wealth penalty, which we name the ``leverage penalty''.

\begin{definition}
    We define the \emph{instantaneuous leverage ratio} before initiating the bet $B_t$ as $\rho_t =  (W_{t-1} + \beta_t)/ W_{t-1}$, which is $\cF_{t-1}$ measurable. We define the \emph{total leverage penalty} up to time $t$ as their product $P_t = \rho_1 \dots \rho_t$. Finally, we call $W_t / P_t $  the \emph{leverage-penalized gross wealth} at time $t$.
Here, while computing $\rho_t$ and $P_t$, we stipulate that $\frac{0}{0}:= 1$ and $0 \cdot \infty = \infty$. 
\end{definition}

\begin{theorem}\label{prop:lvg-rate-mtg}
      The leverage-penalized gross wealth process $\{ W_t/P_t \}$
is a nonnegative {supermartingale under $\Ps$} and {martingale under $\Ph$}. Consequently,
\begin{equation}
    \Ps( \exists t,\ W_t \ge P_t/\alpha ) \le \alpha
\end{equation}
for any $\alpha \in (0,1)$.
\end{theorem}
\begin{proof} Since $\rho_1,\dots, \rho_t$ are all $\cF_{t-1}$-measurable,
    \begin{align}
     \Exph\left(   \frac{W_t}{\rho_1 \dots \rho_t}  \middle| \cF_{t-1} \right ) = \frac{ \Exph( W_t | \cF_{t-1} )}{\rho_1 \dots \rho_t}  
   \stackrel{(*)} =  \frac{W_{t-1} + \beta_t}{\rho_1 \dots \rho_t} =  \frac{W_{t-1}\rho_t}{\rho_1 \dots \rho_t}  =   \frac{W_{t-1}}{\rho_1 \dots \rho_{t-1}},
\end{align}
{and the equality $(*)$ becomes $\le$ if under $\Exps$.} The probability bound follows from \eqref{eqn:time-unif-nb}.
\end{proof}


Importantly, \cref{prop:lvg-rate-mtg} does not require any assumptions other than the set-up \cref{def:nsm-bet2}.
It is worth discussing the effect of the stipulations $\frac{0}{0}:= 1$ and $0 \cdot \infty = \infty$: while we pose no assumptions on the betting procedure, the leverage-penalized process $\{W_t/P_t\}$ does drop to 0 and stays 0 forever once Statistician empties the gross wealth $W_t$, either by betting recklessly ($B_t = 0$, ``betting all-in and losing all''), or by paying back debts with all existing gross wealth leaving zero future capital to bet with ($\beta_t = - W_{t-1}$). The previous method \cref{thm:g-crossing}, in contrast, allows $B_t = 0$ while $\{W_t\}$ crossing a nontrivial boundary $\{g(L_t)\}$.

With \cref{prop:lvg-rate-mtg}, we achive \eqref{eqn:dep-objective} by setting
\begin{equation}
    h(\alpha, w_0,\ell_1,\dots, w_{t-1},\ell_t ) = \frac{1}{\alpha} \cdot \prod_{i=1}^t \frac{w_{i-1}+ \ell_i  -\ell_{i-1} }{w_{i-1}}, \quad (\ell_0 := 0)
\end{equation}
which satisfies $h(\alpha, \cdots) = 1/\alpha$ if all $\ell_i$'s are 0, and $h$ is always finite. That is, the threshold needs not increase when Statistician is not borrowing.

A closer look at $W_t/P_t$, however, reveals that $W_t/P_t$ is not simply penalizing the accrued wealth $W_t$ by the leverage ratios $\rho_1, \dots,  \rho_t$, but itself a wealth supermartingale: {the one that corresponds to the betting \emph{without} borrowing strategy}
\begin{equation}
   { W_t/P_t = B_1 \dots B_t,}
\end{equation}
that is, $\{W_t/P_t\}$ is what the wealth \emph{would have been} if Statistician never borrowed, but used the same {bets $\{ B_t \}$ as they actually do (e.g.\ same betting strategy $\{ \lambda_t \}$ in fair-coin betting)}. {The evolution of these two processes can be compared as follows.}
\begin{align}
     W_0 &= 1,  &  W_t &= (W_{t-1} + \beta_t) \cdot B_t; \label{eqn:in-fact}
     \\
       W_0/P_0 &=  1, &   W_t/P_t &= (W_{t-1}/P_{t-1} + 0) \cdot B_t. \label{eqn:equivalent-non-borrowing-game}
\end{align}

\section{Further generalization I: bargaining and interest rates}
\label{sec:rates}

With the introduction of borrowing into the testing-by-betting game, from \cref{def:nsm-bet} to \cref{def:nsm-bet2}, we essentially drop the assumption that the nonnegative wealth process $\{ W_t \}$ is a supermartingale, by adding the borrowing terms $\{ \beta_t \}$ to the wealth evolution equation $W_t = W_{t-1}\cdot B_t$. As remedies, the two methods explored in the preceding sections both involve first ``reconstructing'' an underlying supermartingale: the net wealth $\{ N_t \}$ in the first method (\cref{sec:gL}), or the leverage ratio-penalized wealth $\{ W_t/P_t \}$ in the second method (\cref{sec:current}). These supermartingales are constructible thanks to the (sub)fairness of the bets $B_t$ under $\Ps$.

In this section, we stretch this line of reasoning even further by adding another term to the wealth evolution equation $W_t = (W_{t-1} + \beta_t)\cdot B_t$ in \cref{def:nsm-bet2}, leading to the following generalization.

\begin{definition}[Borrowed and Bargained NSM Betting]\label[definition]{def:nsm-bet3}
    \normalfont
    Let $(\Omega, \cA, \Pr)$ be a probability space with filtration $\{\cF_t\}_{t \ge 0}$. Let $\{B_t \}_{t \ge 1}$, $\{b_t\}_{t \ge 1}$ be nonnegative processes adapted and predictable w.r.t.\ $\{ \cF_t \}_{t \ge 0}$ respectively,   and $\{ \beta_t \}_{t \ge 1}$ a process predictable w.r.t.\ $\{ \cF_t \}_{t \ge 0}$ such that $\beta_t \ge - W_{t-1}$ for each $t$. Statistician's wealth process $\{ W_t \}_{t \ge 0}$ is defined as $W_0=1$ and
\begin{equation}
     W_t = (W_{t-1} + \beta_t) \cdot  (1+b_t) B_t.
\end{equation}
The supermartingale null $\Ps$ under which $\Exps(B_t | \cF_{t-1} ) \le 1$ for all $t \ge 1$, and the martingale null $\Ph$ under which $\Exph(B_t | \cF_{t-1} ) = 1$ for all $t \ge 1$ are as before.
\end{definition}
Apparently, \cref{def:nsm-bet2} is the special case of \cref{def:nsm-bet3} with $b_t = 0$. With $b_t > 0$, this setup can be interpreted as
 the scenario where the bet takes place alongside a risk-free instrument; or equivalently, the scenario where Statistician has successfully ``bargained'' for enjoying the same potential payoff by risking less capital. 
 
 Let us elucidate this by recalling the fair-coin game example mentioned in \cref{sec:noborrow}, where $X_t = \pm 1$ is the coin-toss outcome and $B_t = 1 + \lambda_t X_t$ is the payoff rate at time $t$. In the previous setup, Casino offers \emph{double-or-nothing} bets on both heads and tails, in which case the bets are fairly priced to the effect that \emph{no arbitrage opportunities exist} for Statistician to win risk-free money. \cref{def:nsm-bet3}, on the other hand, specifies a payoff rate of $(1 + \lambda_t X_t)(1+b_t)$, that is, the bet offered on $X_t$ is \emph{$(2+2b_t)$-times-or-nothing}: for every dollar placed on \emph{either heads or tails},
Casino will pay back $2+2b_t$ dollars if Statistician correctly guessed the side. In this case, we note, Statistician is able to arbitrage by putting $1/2$ on heads and $1/2$ on tails, increasing wealth from 1 to $1+b_t$ regardless of the actual distribution $\Pr$ and regardless of the actual outcome of the coin toss $X_t$. {With \cref{def:nsm-bet3}, Statistician's strategy at time $t-1$ is as follows:} Statistician allocates $|\lambda_t|$ portion of their {post-borrow wealth $W_{t-1} + \beta_t$} into risky bets on heads ($\lambda_t > 0$) or tails ($\lambda_t < 0$), paying off
\begin{equation}
    (W_{t-1} + \beta_t) \cdot  |\lambda_t| (2 + 2 b_t) \id_{ \{ \lambda_t X_t > 0   \} } 
\end{equation}
as per Casino's protocol; and $1-|\lambda_t|$ of their {post-borrow wealth $W_{t-1} + \beta_t$} on the risk-free 50-50 combination of bets on heads and tails, paying off
\begin{equation}
    (W_{t-1} + \beta_t) \cdot  (1-|\lambda_t|)(1+b_t)
\end{equation}
as it does not make sense for Statistician to keep this idle cash uninvested. In total, therefore, the wealth after the payoff is
\begin{equation}
     (W_{t-1} + \beta_t) \cdot   (1 + b_t) \{ 2|\lambda_t| \id_{ \{ \lambda_t X_t > 0   \} }  + (1-|\lambda_t|)\} =  (W_{t-1} + \beta_t) \cdot (1+b_t)(1+\lambda_t X_t),
\end{equation}
{agreeing with \cref{def:nsm-bet3}.

The general setup with bet process $\{B_t\}$ admits a simple economic interpretation: \emph{cash earns interest at rate $b_t$}. Here, $B_t$ is no longer what 1 unit of wealth invested in Statistician's portfolio at time $t-1$ becomes at time $t$; but the $(t-1)$-\emph{present value} of what 1 unit of wealth invested in the portfolio becomes at time $t$, discounted by $1+b_t$.
Statistician's wealth thus becomes $W_{t-1} \cdot B_t \cdot(1+b_t)$ at time $t$, coinciding with \cref{def:nsm-bet3}.
}


The existence of a risk-free instrument leads to the intuition of setting the interest of borrowing at the same rate as the risk-free return rate to prevent pathological possibilities in the system. This idea is mathematically formalized via the first defining the \emph{compound interest liabilities process} to generalize the original definition of liabilities in \eqref{eqn:liabilities},
\begin{equation}
  L_0  = 0, \quad  L_t = (1+b_t)(L_{t-1} + \beta_t),
\end{equation}
and the net wealth as $N_t = W_t - L_t$. Then, we have the following generalization of \cref{prop:doob-wealth} that involves the \emph{present values}  as discounted by the compound interest rate $(1+b_1)\dots(1+b_t)$ of these processes.

\begin{proposition}\label[proposition]{prop:doob-interest}
    {For any $t \ge 1$, it holds that} 
    \begin{equation}\label{eqn:doob-interest}
        \underbrace{\frac{W_t}{(1+b_1)\dots(1+b_t)}}_{W_t'} =  \underbrace{\frac{N_t}{(1+b_1)\dots(1+b_t)}}_{N_t'} +    \underbrace{\frac{L_t}{(1+b_1)\dots(1+b_t)}}_{L_t'}
    \end{equation}
    {where $\{ N_t' \}$ is a supermartingale under $\Ps$ and $\{ L_t' \}$ is predictable}. Further, under $\Ph$,  $\{ N_t' \}$ is a martingale and \eqref{eqn:doob-interest} is the Doob decomposition of the process $\{ W_t' \}$. Additionally, $\{W_t'\}$ evolves as the wealth in a game where the risk-free instrument does not exist:
    \begin{equation}\label{eqn:interest-evolve}
        W_t' =  (W_{t-1}' + \beta_t') \cdot B_t
    \end{equation}
    where $\beta_t ' =L_t' - L_{t-1}'$.
\end{proposition}
\begin{proof}
First, we have
\begin{align}
        & L_t' - L_{t-1}'  = \frac{(1+b_t)(L_{t-1} + \beta_t) -(1+b_t)L_{t-1} }{(1+b_1)\dots(1+b_{t})} = \frac{ \beta_t(1+b_t) }{(1+b_1)\dots(1+b_{t})}
        \\
        = &  \frac{(W_{t-1}+\beta_t)(1+b_t) - (1+b_t)W_{t-1}}{(1+b_1)\dots(1+b_{t})}  = W_t'/B_t - W_{t-1}',
    \end{align}
    which implies  \eqref{eqn:interest-evolve}. This further implies that, under $\Ph$, $L_t' - L_{t-1}' = \Exph(W_t' | \cF_{t-1}) - W_{t-1}'$,
    concluding that $N_t' + L_t'$ is the Doob decomposition of $W_t'$. The case under $\Ps$ is analogous.
\end{proof}

This decomposition in \cref{prop:doob-interest} immediately enables the generalization of \cref{thm:g-crossing,prop:lvg-rate-mtg} in this setting.

\begin{theorem} Denote $R_t = \prod_{i=1}^t(1+b_i)$. Then:
    \begin{enumerate}
        \item Under the same assumptions as \cref{thm:g-crossing},
        \begin{equation}\label{eqn:g-bound-2}
    \Ps( \exists t, \ W_t \ge R_t \cdot g(L_t / R_t) ) \le 1 - \frac{g(0) - 1}{g(0)}\exp \left( -\int_0^\infty \frac{1}{g(\ell)} \d \ell \right);
\end{equation}
\item With $P_t = \prod_{i=1}^t \frac{W_{i-1}' + \beta_i'}{W_{i-1}'} = \prod_{i=1}^t \frac{W_{i-1} + \beta_i}{W_{i-1}}$,
\begin{equation}
    \Ps( \exists t,\ W_t \ge R_t P_t/\alpha ) \le \alpha
\end{equation}
for any $\alpha \in (0,1)$.
    \end{enumerate}
\end{theorem}
\begin{proof}
    The only non-straightforward step is that $ \frac{W_{t-1}' + \beta_t'}{W_{t-1}'} =  \frac{W_{t-1} + \beta_t}{W_{t-1}}$ for any $t$. To see this,
    \begin{equation}
        \beta_t' R_{t-1} = \frac{L_t}{1+b_t} - L_{t-1} = \beta_t. \qedhere
    \end{equation}
\end{proof}

It is worth noting that if Statistician refrains from borrowing in this case, i.e., $P_t = 1$, there is then a clear analogy between the second bound above $\Ps( \exists t,\ W_t \ge R_t /\alpha ) \le \alpha$ and the previous bound $  \Ps( \exists t,\ W_t \ge P_t/\alpha ) \le \alpha$ in \cref{prop:lvg-rate-mtg}. This suggests the following (informal) way to relate borrowing and bargaining: they both allow ``free cash'', at the growth rate of $P_t$ and $R_t$ respectively. For the first bound above, if Statistician does not borrow, we observe the same effect (as \cref{thm:g-crossing}) of a loosened Ville's inequality on the NSM $\{ W_t/R_t \}$, as we noted in \cref{prop:g-deficiency}.

\section{Further generalization II: varying the initial conditions}
\label{sec:initial_conditions}

Throughout the preceding sections, we have assumed a baseline where Statistician begins with initial gross wealth $W_0 = 1$ and initial liabilities $L_0 = 0$. In this section, we examine the consequences of relaxing these assumptions across the problem space 
$$\{\text{\cref{thm:g-crossing}}, \text{\cref{prop:lvg-rate-mtg}}\} \times \{W_0 \in \mathbb{R}\} \times \{L_0 \in \mathbb{R}\}.$$

We can immediately prune several regions of this space. First, we do not consider $W_0 < 0$, as this violates the meaning of ``gross'' wealth and arguably defies the purpose of introducing the distinction
between the supermartingale net wealth $N_t$ and the nonnegative gross wealth $W_t$. Second, an initial debt $L_0 = \ell > 0$ is a trivial translation. Because $L_t = L_0 + \sum_{i=1}^t \beta_i$, one can equivalently define a new liabilities process $L_t' = L_t - L_0$ starting at $0$ and shift the threshold curves accordingly without altering the structural form of any bounds. Third, scaling $W_0$ to an arbitrary strictly positive value $w > 0$ simply rescales the betting processes natively.

The case where Statistician arrives at Casino with strictly zero capital ($W_0 = 0, L_0 = 0$), however, is particularly interesting, and it fundamentally distinguishes the two methods. The leverage penalty framework (\cref{prop:lvg-rate-mtg}) fails entirely. The initial leverage ratio becomes $\rho_1 = (W_0 + \beta_1)/W_0 = \infty$, rendering the cumulative penalty infinite and any finite rejection threshold unreachable. This failure is philosophically consistent: \cref{prop:lvg-rate-mtg} adjusts the realized gross wealth by the counterfactual borrow-free wealth, making $W_0 = 0$ equivalent to a bettor starting with zero capital and barred from borrowing, whose wealth is permanently at zero.

In contrast, the liabilities curve framework (\cref{thm:g-crossing}) naturally accommodates $W_0 = 0$. We formalize this intriguing parallel universe below.

\subsection{The zero-wealth universe ($W_0 = 0$)}

We define the ``$\$0$-universe'' as the scenario where $W_0 = 0$ and Statistician, now referre to as the ``$\$0$-Statistician'' to differentiate from the previous ``$\$1$-Statistician'', must borrow ($\beta_1 > 0$) to place their initial bet. 

\begin{theorem}[Zero-wealth liabilities curve bound] \label{thm:zero-wealth-bound}
    Assume $W_0=0$, $L_0=0$, and that the liabilities process $\{L_t\}$ is non-decreasing. For any function $g_0:[0,\infty) \to [0,\infty)$ such that $\ell \mapsto g_0(\ell) - \ell$ is increasing, the probability of falsely rejecting the supermartingale null $\Ps$ is bounded by:
    \begin{equation}\label{eqn:bound-0-univ}
        \Ps \left( \exists t,\ W_t \ge g_0(L_t) \right) \le 1 - \exp\left(-\int_0^\infty \frac{1}{g_0(x)} \d x \right).
    \end{equation}
\end{theorem}
The proof follows identically from substituting $W_0=0$ into the initial auxiliary supermartingale $K_0$ in the proof of \cref{thm:g-crossing}. Comparing \eqref{eqn:bound-0-univ} to the standard ``$\$1$-universe'' ($W_0=1$) bound \eqref{eqn:g-bound}:
 \begin{equation}\label{eqn:g1-bound}
    \Ps( \exists t, \ W_t \ge g_1(L_t) ) \le 1 - \frac{g_1(0) - 1}{g_1(0)}\exp \left( -\int_0^\infty \frac{1}{g_1(x)} \d x \right),
\end{equation}
we make the following observations.

First, the $\$0$-universe bound \eqref{eqn:bound-0-univ} removes the ``intercept penalty'' term $\frac{g_1(0) - 1}{g_1(0)}$ in the $\$1$-universe bound \eqref{eqn:g1-bound}, and subsequently removes the ``deficiency'' described in \cref{prop:g-deficiency} (but we shall soon see the $\$0$-incarnation of the same limit in \cref{sec:g0-no-free-lunch}).
To calibrate both universes to a significance level $\alpha$, the $\$1$-universe bound mandates that the threshold curve $g_1(x)$ strictly satisfies $g_1(0) \ge 1/\alpha$. For example, at $\alpha = 0.05$, the $\$1$-Statistician must employ a curve that starts at $g_1(0) \ge 20$. The $\$0$-universe bound circumvents this entirely; it places no lower bound constraint on $g_0(0)$. The $\$0$-Statistician could validly employ a curve starting at an arbitrarily small $g_0(0) > 0$, provided the integral converges to $-\log(1-\alpha)$. 

Second, we note the following two-way emulation argument between the $\$0$- and the $\$1$-universes: the two Statisticians may artificially adjust their initial borrowing and invoke the bound from the other universe, arriving at the same curve-crossing probability.

 \begin{itemize}
        \item \textbf{$\$1$-Statistician emulating the $\$0$-universe:} Suppose the $\$1$-Statistician (having picked a $g_1$) artificially relinquishes their initial $\$1$, instead borrowing it at $t=1$ ($\beta_1^{(0)} = \beta_1^{(1)} + 1$). Their emulated liabilities shift by $+1$. To trigger identical rejection paths, their original curve $g_1$ maps to $g_0(x) = g_1(x - 1)$ for $x \ge 1$. To apply the $\$0$-universe theorem, $g_0$ must be extended to $[0,1)$ such that $g_0(x) - x$ is increasing. Since $g_0$ on $[0,1)$ is unvisited, we may pick the largest choice
      $g_0(x) = x + g_1(0) - 1$ on $x \in [0,1)$ to obtain the tightest bound (i.e.\ smallest curve-crossing probability). The extra integral over this unvisited interval corresponds exactly to the ``intercept penalty'' mentioned above
        \[
      \int_0^1 \frac{1}{g_0(x)} \d x =  \int_0^1 \frac{1}{x+g_1(0) - 1} \d x = -\log\left(\frac{g_1(0) - 1}{g_1(0)}\right).
        \]
        Substituting this into the $\$0$-universe bound exactly recovers \emph{both sides} of the $\$1$-universe bound (\cref{thm:g-crossing}):
        \[
        1 - \exp\left(-\int_0^\infty \frac{1}{g_0(x)} \d x \right) = 1 - \frac{g_1(0) - 1}{g_1(0)}\exp \left( -\int_0^\infty \frac{1}{g_1(x)} \d x \right).
        \]
        Thus, we may say that the intercept penalty $\frac{g_1(0)-1}{g_1(0)}$ is incurred due to the debt trajectory $L_t \in [0,1]$ in the $\$0$-universe \emph{``skipped'' in the emulation}.
        \item \textbf{$\$0$-Statistician emulating the $\$1$-universe:} Suppose the $0$-Statistician (having picked a $g_0$) artificially gifts themselves $\$1$ in initial capital, reducing their emulated borrowing by $\$1$ ($\beta_1^{(1)} = \beta_1^{(0)} - 1$). Their emulated liabilities shift by $-1$. To trigger identical rejection paths, their original curve $g_0$ maps to $g_1(x) = g_0(x + 1)$. Applying the $\$1$-universe theorem evaluates their bound as $1 - \frac{g_0(1) - 1}{g_0(1)} \exp\left(-\int_1^\infty \frac{1}{g_0(y)} \d y\right)$. By the constraint that their original $g_0(x) - x$ is increasing, it holds that $g_0(x) \le x + g_0(1) - 1$ on $[0,1]$. Thus, the true integral $\int_0^1 \frac{1}{g_0(x)} \d x \ge -\log\left(\frac{g_0(1) - 1}{g_0(1)}\right) = -\log\left(\frac{g_1(0) - 1}{g_1(0)}\right)$. Therefore,
        \[
        1 - \exp\left(-\int_0^\infty \frac{1}{g_0(x)} \d x \right) \ge 1 - \frac{g_1(0) - 1}{g_1(0)}\exp \left( -\int_0^\infty \frac{1}{g_1(x)} \d x \right),
        \]
        meaning that the bound they compute via $\$1$-universe emulation implies the $\$0$-universe bound they could have achieved, also recovering both sides of the bound.
    \end{itemize}

\begin{table}[htbp]
    \centering
    \begin{tabular}{p{0.3\textwidth} p{0.3\textwidth} p{0.3\textwidth}}
        \toprule
        \textbf{Concept} & \textbf{$\$1$-Statistician emulating $\$0$} & \textbf{$\$0$-Statistician emulating $\$1$} \\
        \midrule
        Initial capital adjustment & $W_0^{(0)} = W_0^{(1)} - 1$ & $W_0^{(1)} = W_0^{(0)} + 1$ \\
        First borrowing shift & $\beta_1^{(0)} = \beta_1^{(1)} + 1$ & $\beta_1^{(1)} = \beta_1^{(0)} - 1$ \\
        Wealth shift ($t \ge 1$) & None: $W_t^{(0)} = W_t^{(1)}$ & None: $W_t^{(0)} = W_t^{(1)}$ \\
        Liabilities shift ($t \ge 1$) & $L_t^{(0)} = L_t^{(1)} + 1$ & $L_t^{(1)} = L_t^{(0)} - 1$ \\
        Threshold curve mapping & $g_0(x) = g_1(x-1)$ for $x\ge 1$, linear extrapolation for $x \in [0,1)$ & $g_1(x) = g_0(x+1)$ for $x \ge 0$ \\
        \bottomrule
    \end{tabular}
    \caption{Correspondences in the two-way emulation argument between the $\$0$- and $\$1$-universes.}
    \label{tab:two-way-emulation}
\end{table}
    
We summarize the two-way emulation argument in \cref{tab:two-way-emulation}. In this argument, the procedure ``$\$i$-Statistician emulating the $\$(1-i)$-universe'' shows that the $\$(1-i)$-universe bound implies the $\$i$-universe bound. When $i=1$, the implication does not require any assumption on the $\$1$-universe's process. When $i=0$, the implication requires that $\beta_1^{(0)}\ge 1$. In this sense, we might say that the $\$0$-universe bound \cref{thm:zero-wealth-bound} is slightly stronger than the $\$1$-universe bound \cref{thm:g-crossing}: \emph{any instance of the $\$1$-universe bound can be derived from an instance of the $\$0$-universe bound; any instance of the $\$0$-universe bound with initial borrowing $\ge 1$ can be derived from an instance of the $\$1$-universe bound}. Once again, however, we have argued that all $\$w$-universes ($w > 0$) are equivalent, and one may argue that ``any instance of the $\$0$-universe bound with initial borrowing $\ge w$ can be derived from an instance of the $\$w$-universe bound'', where $w$ may be arbitrarily small. However, fixing \emph{a} $\$w$-universe there are always some $\$0$-Statisticians unemulable in this universe.
The $\$0$-universe thus provides a scale-agnostic treatment of the testing-by-betting thesis.

The discussions above on the $\$0$-universe versus the $\$1$-universe apply \emph{mutatis mutandis} to the bargaining setting as well. We omit the corresponding statements from this paper.

Finally, some quadratic examples of $g_1$ functions are provided in \cref{sec:g-inventory}, together with the illustration of ``shifting'' in the emulation argument.
The rest of this section contains a few novelties around the $\$0$-universe.

\subsection{No free lunch in the $\$0$-universe}\label{sec:g0-no-free-lunch}

The equivalences established above raise an intuitive question regarding the \emph{timing} of borrowing. If the $\$1$-Statistician starts with $\$1$ and borrows nothing ($\beta_1^{(1)} = 0$), they face a massive baseline threshold $g_1(0) \ge 1/\alpha$. If the $\$0$-Statistician instead borrows that $\$1$ \emph{gradually}---for instance, borrowing small increments $\delta$ over many steps until their total debt reaches $L_N = 1$---can they exploit the $\$0$-universe's arbitrarily low starting threshold $g_0(0)$ to reject the null hypothesis with less evidence?

 Remarkably, the answer is no; the mathematical constraints of the $\$0$-universe perfectly preserve the $1/\alpha$ evidence barrier with the following no-free-lunch property.

\begin{proposition}\label[proposition]{prop:g0-bound-nfl}
    Under the $\$0$-universe bound (\cref{thm:zero-wealth-bound}), any valid threshold curve $g_0(x)$ calibrated to a significance level $\alpha \in (0,1)$ must satisfy
    \begin{equation}
        g_0(w) \ge \frac{w}{\alpha}
    \end{equation}
    for all $w > 0$. In particular, $g_0(1) \ge 1/\alpha$.
\end{proposition}

\begin{proof}
    By \cref{thm:zero-wealth-bound}, valid threshold curves must satisfy $\int_0^\infty \frac{1}{g_0(x)} \d x \le -\log(1-\alpha)$. Furthermore, the constraint that $x \mapsto g_0(x) - x$ is increasing implies that $g_0(x) \le g_0(w) - (w - x)$ for $x \in [0, w]$. Therefore,
    \begin{equation}
    \log\left( \frac{g_0(w)}{g_0(w) - w} \right) =  \int_0^w \frac{1}{g_0(w) - w + x} \d x  \le \int_0^\infty \frac{1}{g_0(x)} \d x \le -\log(1-\alpha).
    \end{equation}
    This implies that
    \begin{equation}
        \frac{g_0(w)}{g_0(w) - w} \le \frac{1}{1-\alpha} \implies g_0(w) \ge \frac{w}{\alpha}. \qedhere
    \end{equation}
\end{proof}

To conclude, the $\$1$-Statistician evaluates their initial $\$1$ of capital against the threshold $g_1(0) \ge 1/\alpha$; whereas the $\$0$-Statistician, attempting to ``sneak under'' the barrier by borrowing $\$1$ gradually, evaluates their wealth against $g_0(1)$. Yet, the steepness constraint forces the curve $g_0$ to aggressively spike precisely proportional to the integral cost of traversing the debt interval $[0, 1]$, universally locking the minimum required evidence at $1/\alpha$.

\subsection{$L_0 = 1$ in the $\$1$-universe}

While we have argued before that $L_0$ does not matter, it is helpful now to consider viewing the initial $\$1$ in the $\$1$-universe not as ``gifted'' ($L_0=0$), but as already borrowed ($L_0=1, W_0=1$), as this adds another perspective (although this is logically redundant) to the two-way emulation picture.

We term this the ``borrowed $\$1$-universe''. In this setting, the two-way emulation with the $\$0$-universe becomes even simpler because the liabilities trajectories match exactly for $t \ge 1$: if the $\$0$-Statistician borrows $\$1$ at $t=1$ to match the $\$1$-Statistician's wealth, both Statisticians possess $L_1 = 1$. Consequently, there is no spatial shift in the threshold curves; the equivalent threshold functions in the two universes satisfy $g_0(x) = g_1(x)$, only differing in that $g_i$ is defined on $[i,\infty)$.
    
    The $\$1$-universe bound for $L_0=1, W_0=1$ reads
    \begin{equation}\label{eqn:g1b-bound}
   \Ps( \exists t, \ W_t \ge g_1(L_t) ) \le 1 - \frac{g_1(1) - 1}{g_1(1)}\exp \left( -\int_1^\infty \frac{1}{g_1(x)} \d x \right),
\end{equation}
   \begin{itemize}
        \item \textbf{Borrowed $\$1$ emulating $\$0$:} The borrowed-$\$1$ Statistician pretends they started with $\$0$ and borrowed their initial $\$1$ at $t=1$. To apply the $\$0$-universe bound ($1-\exp(-\int_0^\infty \dots)$), their original curve $g_1(x)$ (defined for $x \ge 1$) must be extended over the unvisited debt interval $[0,1)$. The optimal extension satisfying the steepness constraint is $g_0(x) = x + g_1(1) - 1$. The integral cost over this unvisited interval is exactly $\int_0^1 \frac{1}{x + g_1(1) - 1} \d x = \log(\frac{g_1(1)}{g_1(1)-1})$. Evaluating the $\$0$-universe bound then yields
         \[
        1 - \exp\left(-\int_0^\infty \frac{1}{g_0(x)} \d x \right) = 1 - \frac{g_1(1) - 1}{g_1(1)}\exp \left( -\int_1^\infty \frac{1}{g_1(x)} \d x \right),
        \]
        which exactly recovers their true bound \eqref{eqn:g1b-bound} in both sides.
        \item \textbf{$\$0$ emulating borrowed $\$1$:} The $\$0$-Statistician borrows $\$1$ immediately and pretends this was their initial state ($L_0=1, W_0=1$). Applying the borrowed-$\$1$ universe bound evaluates to $1 - \frac{g_0(1) - 1}{g_0(1)}\exp(-\int_1^\infty \frac{1}{g_0(x)} \d x)$. By the steepness constraint, $\int_0^1 \frac{1}{g_0(x)} \d x \ge \log(\frac{g_0(1)}{g_0(1)-1})=\log(\frac{g_1(1)}{g_1(1)-1})$, which implies \[
        1 - \exp\left(-\int_0^\infty \frac{1}{g_0(x)} \d x \right) \ge 1 - \frac{g_1(1) - 1}{g_1(1)}\exp \left( -\int_1^\infty \frac{1}{g_1(x)} \d x \right),
        \]
         meaning that the bound they compute via borrowed $\$1$-universe emulation implies the $\$0$-universe bound they could have achieved, also recovering both sides of the bound.
    \end{itemize}
    This confirms that the two universes are perfectly mathematically equivalent. See \cref{tab:two-way-emulation-2}.

    \begin{table}[htbp]
    \centering
    \begin{tabular}{p{0.3\textwidth} p{0.3\textwidth} p{0.3\textwidth}}
        \toprule
        \textbf{Concept} & \textbf{Borrowed $\$1$-Statistician emulating $\$0$} & \textbf{$\$0$-Statistician emulating borrowed $\$1$} \\
        \midrule
        Initial capital adjustment & $W_0^{(0)} = W_0^{(1)} - 1$ & $W_0^{(1)} = W_0^{(0)} + 1$ \\
        Initial liabilities adjustment & $L_0^{(0)} = L_0^{(1)} - 1$ & $L_0^{(1)} = L_0^{(0)} + 1$ \\
        First borrowing shift & $\beta_1^{(0)} = \beta_1^{(1)} + 1$ & $\beta_1^{(1)} = \beta_1^{(0)}  - 1$ \\
        Wealth shift ($t \ge 1$) & None: $W_t^{(0)} = W_t^{(1)}$ & None: $W_t^{(0)} = W_t^{(1)}$ \\
        Liabilities shift ($t \ge 1$) & None: $L_t^{(0)} = L_t^{(1)}$ & None:   $L_t^{(1)} = L_t^{(0)}$ \\
        Threshold curve mapping & $g_0(x) = g_1(x)$ for $x \ge 1$, linear extrapolation for $x \in [0,1)$ & $g_1(x) = g_0(x)$ for $x \ge 1$ \\
        \bottomrule
    \end{tabular}
    \caption{Correspondences in the two-way emulation argument between the $\$0$-universe and borrowed $\$1$-universe.}
    \label{tab:two-way-emulation-2}
\end{table}

\section{Summary}
We study in this paper the extension of ``testing by betting", a central topic in game-theoretic statistics, into the setting where borrowing is allowed. 
We derived two types of sequential rejection rules that generalize the classical ``reject when wealth above $1/\alpha$'' rule, with different consequences.

 In the first type of rules, Statistician adjusts the wealth made from betting by the total borrowed amount, an idea mathematically formalized by a curve-crossing bound on the probability of the gross wealth ever crossing a fixed function of the liabilities:
 \begin{equation}
     \Ps( \exists t, W_t \ge g(L_t) ) \le 1 - \frac{g(0)-1}{g(0)} e^{-\int_0^\infty \frac{1}{g(x)} \d x }.
 \end{equation}
This incurs \emph{two} prices: even if Statistician does not borrow, the \emph{optionality} of borrowing increases the rejection threshold ($g(0) > 1/\alpha$); if Statistician does borrow, the threshold grows much faster than the borrowed amount does ($g(x) \gtrsim x$). The benefit, we note, is that it allows borrowing \emph{after bankruptcy} to bet again, without ruining the curve-crossing possibility and hence power.

 In the second type of rules, Statistician adjusts the wealth by $P_t$, the product of past leverage ratios:
 \begin{equation}
     \Ps( \exists t, W_t/P_t \ge 1/\alpha ) \le \alpha.
 \end{equation}
 This adjusted process $\{ W_t/P_t \}$ is equivalent to the e-process of betting without borrowing $W_t/P_t = B_1\dots B_t$. Thus, when not borrowing ($P_t = 1$), this approach recovers the ``$\Ps( \exists t, W_t \ge 1/\alpha ) \le \alpha $'' no-borrow bound.
 Here, the optionality of borrowing comes at no cost or gain in comparison to some no-borrow test instance. However, this approach does not allow borrowing after bankruptcy (while retaining any statistical power), which forces $W_t/P_t = 0$ forever.
 
This dichotomy demonstrates an intriguing trade-off in our generalized testing-by-betting framework.

\subsection*{Acknowledgements} 

The authors acknowledge Muriel F.\ P\'erez-Ortiz for helpful discussions.

\bibliography{main}

\newpage

\appendix

\section{On the $g(L_t)$-crossing bound (\cref{thm:g-crossing})}

\subsection{Monotonicity of borrowing}\label[appendix]{sec:why-not-move-back}

In this section, we argue that we shall not expect \cref{thm:g-crossing} to hold when $\{L_t\}$ is allowed to decrease. Recall \cref{thm:g-crossing} states that:

\begin{theorem}\label{thm:g-crossing-2} Assume the liabilities process $\{L_t\}$ is non-decreasing.
Then, for any  function $g:[0,\infty) \to [1,\infty)$ such that $\ell \mapsto g(\ell) - \ell$ is increasing,
    \begin{equation}
    \Ps( \exists t, \ W_t \ge g(L_t) ) \le 1 - \frac{g(0) - 1}{g(0)}\exp \left( -\int_0^\infty \frac{1}{g(\ell)} \d \ell \right).
\end{equation}
\end{theorem}

If $L_t$ is allowed to decrease, Statistician can always pay back all debt and bet 0 dollars on the next observation (i.e.\ $B_{t+1}=1$), resulting
\begin{equation}
    W_t \longrightarrow W_{t+1} = W_t - L_t,\quad L_t \longrightarrow L_{t+1} = 0.
\end{equation}
We note that
\begin{equation}
  \boxed{ \{ W_{t+1} \ge g(0) \} = \{ W_{t+1} \ge g(L_{t+1}) \} \supsetneq \{  W_t \ge g(L_t) \}. }
\end{equation}
This is because, when $W_t \ge g(L_t)$ holds, due to the monotonicity of $\ell \mapsto g(\ell) - \ell$:
\begin{equation}
   W_t \ge g(L_t) = g(L_t) - L_t + L_t \ge g(0) - 0 + L_t = g(0) + L_t,
\end{equation}
hence $W_{t+1} = W_t - L_t \ge g(0)$.

The boxed inclusion states that, when
\begin{equation}
    \omega \notin \{  W_t \ge g(L_t) \}, \quad \omega \in \{ W_{t+1} \ge g(0) \},
\end{equation}
it is possible to \emph{go from non-rejection to rejection without observing any data} (i.e.\ no randomness involved).

In fact, the monotonicity of $\ell \mapsto g(\ell) - \ell$ ensures that Statistician cannot go from non-rejection to rejection by just borrowing a large sum of money and doing nothing:
\begin{equation}
    \omega \notin \{  W_t \ge g(L_t) \} \implies \omega \notin \{  W_t + \beta_{t+1} \ge g(L_t  + \beta_{t+1}) \}.
\end{equation}

\subsection{An inventory of $g$ functions} \label[appendix]{sec:g-inventory}

In this subsection, we provide some closed-form functions $g$ such that the right hand side of the bound \eqref{eqn:g-bound} is calibrated to a pre-specified $\alpha$. Let us set $\alpha=0.05$ for illustration. Recalling our generalization in \cref{sec:initial_conditions}, we include both $g_1$, the threshold curve in the original $\$1$-universe \eqref{eqn:g-bound}; and $g_0$, the target curve in the $\$0$-universe \eqref{eqn:bound-0-univ}.

\begin{figure}[htbp]
    \centering
    \includegraphics[width=\textwidth]{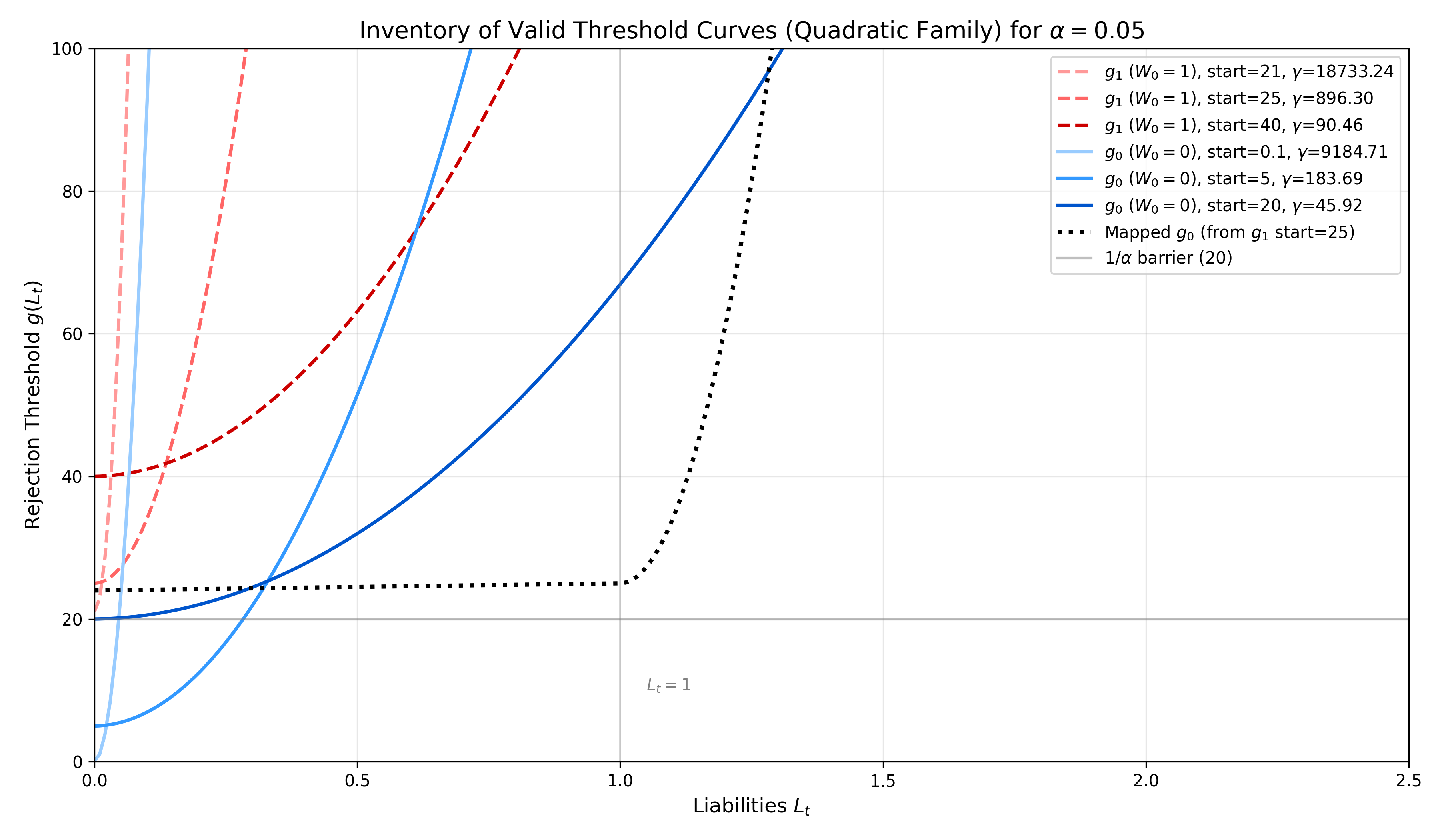}
    \caption{Inventory of valid threshold curves for both universes, calibrated to $\alpha=0.05$.}
    \label{fig:curves-inventory}
\end{figure}

To satisfy these bounds, threshold curves must satisfy the steepness condition ($g_{0,1}'(\ell) \ge 1$) and have an integral of their reciprocal bounded by the targeted confidence limits. We construct an inventory of curves within the \emph{quadratic family}: $g_{0,1}(\ell) = c + \ell + \gamma \ell^2$. For a chosen intercept $c$, the curvature $\gamma$ is numerically calibrated such that the analytical Type-I error bounds perfectly evaluate to $\alpha=0.05$. These curves are in \cref{fig:curves-inventory}.

\paragraph{The $\$1$-universe ($W_0=1$)}
The threshold $g_1(\ell)$ must satisfy $\int_0^\infty \frac{1}{g_1(\ell)} \d\ell \le -\log(1-\alpha) - \log\left(\frac{g_1(0)}{g_1(0)-1}\right)$.
Because of the ``intercept penalty'', any valid $g_1$ must begin strictly above $1/\alpha = 20$.
\begin{itemize}
    \item \textbf{$c=21$}: Requires massive curvature ($\gamma \approx 5.3$) to converge quickly enough, as $21$ is dangerously close to the $20$ barrier.
    \item \textbf{$c=25$}: Requires moderate curvature ($\gamma \approx 0.77$).
    \item \textbf{$c=40$}: Allows for very gentle curvature ($\gamma \approx 0.12$).
\end{itemize}

\paragraph{The $\$0$-universe ($W_0=0$)}
The threshold $g_0(\ell)$ must simply satisfy $\int_0^\infty \frac{1}{g_0(\ell)} \d\ell \le -\log(1-\alpha)$.
Without the intercept penalty, the curve can begin literally anywhere, including near zero.
\begin{itemize}
    \item \textbf{$c=0.1$}: Highly aggressive starting point, needs high curvature ($\gamma \approx 4.0$).
    \item \textbf{$c=5$}: Mid-range starting point, moderate curvature ($\gamma \approx 1.2$).
    \item \textbf{$c=20$}: Begins exactly at the $1/\alpha$ barrier, requiring gentle curvature ($\gamma \approx 0.28$).
\end{itemize}

\paragraph{Two-Way Emulation Mapping} Let us also graphically illustrate the emulation argument in \cref{sec:initial_conditions}.
Notice the dotted black line in \cref{fig:curves-inventory}: this is the mapping $g_0(x) = g_1(x-1)$ with the maximally steep unvisited extension $g_0(x) = x + c - 1$ mapped from the $\$1$-Statistician with $c=25$.
It perfectly constructs a valid $\$0$-universe curve, with the integral over the $[0,1]$ segment exactly matching the mathematical ``intercept penalty'' $\log(25/24)$.

\subsection{Proof of \cref{thm:g-crossing}}\label[appendix]{sec:pf-gcross}

We now prove \cref{thm:g-crossing}. The proof of the inequality \eqref{eqn:g-bound} involves a technique similar to the one developed by \citet[Section~2.2]{KOOLEN2026110577}, in particular the continuous time version, but with time $t$ replaced by liabilities $L_t$; whereas the proof of the tightness \eqref{eqn:g-bound-tight} involves picking a martingale game instance with infinitesimal borrowings $\{\beta_t\}$ where \eqref{eqn:g-bound} holds with arbitrarily small margin.

\begin{proof}[Part I: Proof of the inequality \eqref{eqn:g-bound}] Let $h(\ell) := g(\ell) - \ell$.
 We define the auxiliary process
\[
  K_t ~:=~ 1- \frac{h(L_t) - N_t }{h(L_t) + L_t} s(L_t)
  \qquad
  \text{where}
  \qquad
  s(L) ~:=~ e^{-\int_L^\infty \frac{1}{g(\ell)} d \ell}
  .
\]
Observe that, for any $L, \beta \ge 0$, by the monotonicity of $h$,
\begin{equation}\label{eq:ineq}
  \frac{s(L + \beta)}{s(L)}
  ~=~
  e^{\int_L^{L+\beta} \frac{1}{h(\ell) + \ell} \d \ell}
  ~\ge~
  e^{\int_L^{L+\beta} \frac{1}{h(L+\beta) + \ell} \d \ell}
  ~=~
  1 + \frac{\beta}{h(L+\beta)+L}
  .
\end{equation}
Let us demonstrate that $\{ K_t \}$ is a nonnegative supermartingale \emph{before} crossing 1. That is, letting $\tau$ be the stopping time
\begin{equation}
    \tau = \min\{ t :  M_t \ge g(L_t) \} = \min\{ t :  N_t \ge h(L_t) \} = \min\{t : K_t \ge 1 \},
\end{equation}
the stopped process $\{ K^\tau_t \}$ is a nonnegative supermartingale.

To see that it is nonnegative,
\begin{equation}\label{eqn:NL-comparison}
   h(L_t) - N_t  = h(L_t) + L_t - W_t \le  h(L_t) + L_t
\end{equation}
and
\begin{equation}
 0 <   s(L_t) \le 1.
\end{equation}
To see that it is a supermartingale under $\Ps$, define $\beta_{t+1}^* = \beta_{t+1} \id_{\{ t < \tau \}} $ and $B_{t+1}^* = B_{t+1} \id_{\{ t < \tau \}} + \id_{\{t \ge \tau \}}$. Then $\beta_{t+1}^* = L_{t+1}^\tau - L_t^\tau$ and $W_{t+1}^\tau = (W_t^\tau + \beta_{t+1}^* \id_{\{ t < \tau \}})B_{t+1}^*$. We have
\begin{align*}
  &\Exps\left[K_{t+1}^\tau \middle| \mathcal F_t\right]
  \\
  ~=~&
  \Exps\left[
    1- \frac{h(L_{t+1}^\tau) - N_{t+1}^\tau }{h(L_{t+1}^\tau) + L_{t+1}^\tau} s(L_{t+1}^\tau)
    \middle| \mathcal F_t\right]
  \\
  ~=~&
  \Exps\left[
    1- \frac{h(L_t^\tau + \beta_{t+1}^*) - ((W_t^\tau + \beta_{t+1}^*) B_{t+1}^* - (L_t^\tau + \beta_{t+1}^*))}{h(L_t^\tau + \beta_{t+1}^*) + (L_t^\tau + \beta_{t+1}^*)} s(L_t^\tau + \beta_{t+1}^*)
    \middle| \mathcal F_t\right]
  \\
  ~\stackrel{\clap{\text{\tiny (Def~\ref{def:nsm-bet})}}}{\le}~&
  1- \frac{h(L_t^\tau + \beta_{t+1}^*) - ((W_t^\tau + \beta_{t+1}^*) - (L_t^\tau + \beta_{t+1}^*))}{h(L_t^\tau + \beta_{t+1}^*) + (L_t^\tau + \beta_{t+1}^*)} s(L_t^\tau + \beta_{t+1}^*)
  \\
  ~=~&
  1- \frac{h(L_t^\tau + \beta_{t+1}^*) - N_t^\tau}{h(L_t^\tau + \beta_{t+1}^*) + (L_t^\tau + \beta_{t+1}^*)} s(L_t^\tau + \beta_{t+1}^*)
  \\
  ~\stackrel{\eqref{eq:ineq}}{\le}~&
  1- \frac{h(L_t^\tau + \beta_{t+1}^*) - N_t^\tau}{h(L_t^\tau+\beta_{t+1}^*)+L_t^\tau} s(L_t^\tau)
  \\
  ~\stackrel{\eqref{eqn:NL-comparison}}{\le}~&
  1- \frac{h(L_t^\tau) - N_t^\tau}{h(L_t^\tau)+L_t^\tau} s(L_t^\tau)
  \\
  ~=~&
  K_t^\tau
\end{align*}
Finally, $M_t \ge g(L_t) \iff N_t \ge h(L_t) \iff K_t \ge 1$. So by the standard Ville's inequality \eqref{eqn:time-unif-nb},
\[
  \Ps \left\{
    \exists t : M_t \ge g(L_t)
  \right\}
  ~=~
  \Ps \left\{
    \exists t : K_t \ge 1
  \right\}
   ~=~
  \Ps \left\{
    \exists t : K_t^\tau \ge 1
  \right\}
  ~\stackrel{\text{\eqref{eqn:time-unif-nb}}}{\le}~
  \Exps[K_0]
  ~=~
  1- \frac{g(0) - 1}{g(0)} e^{-\int_0^\infty \frac{1}{g(\ell)} d \ell}
  .
\]
This finishes the inequality \eqref{eqn:g-bound}.
\end{proof}

\begin{proof}[Part II: Proof of the tightness \eqref{eqn:g-bound-tight}]

\underline{Step 1: An anti-concentration inequality.} We generalize Lemma~33 in \cite{ramdas2020admissible}, an anti-concentration counterpart to Ville's inequality for nonnegative martingales, to processes that may be called ``almost submartingales'' (see e.g.\ \cite{robbins1971convergence} for the concept of almost supermartingales).
\begin{lemma} \label[lemma]{lem:almost_sub_anti}
Let $\{K_t\}$ be an adapted nonnegative process with non-random $K_0$. Suppose there exists a predictable sequence $b_t \ge 0$ such that $\Exp[K_t \mid \cF_{t-1}] \ge K_{t-1} - b_t$. Let $Y_t = K_t / K_{t-1}$ (with $0/0=1$). Assume that $K_\infty = 0$ almost surely on the event $\{\sup_t K_t < 1\}$. Assume further that there exists $\varepsilon > 0$ such that for any $t \in \mathbb{N}$ and any $\cF_{t-1}$-measurable random variable $\beta \ge 1$, 
\[
\Exp[Y_t \mid \cF_{t-1}, Y_t \ge \beta] \le \beta (1+\varepsilon).
\]
Then, for $\tau = \inf\{t : K_t \ge 1\}$, we have
\[
\Pr\left( \sup_{t\ge 0} K_t \ge 1 \right) \ge \frac{K_0 - \Exp\left[\sum_{t=1}^\tau b_t\right]}{1+\varepsilon}.
\]
\end{lemma}
\begin{proof}
The process $S_t = K_t + \sum_{i=1}^t b_i$ is a nonnegative submartingale. By the optional stopping theorem, $\Exp[S_{\tau \wedge t}] \ge S_0 = K_0$, which implies $\Exps[K_{\tau \wedge t}] \ge K_0 - \Exp[\sum_{i=1}^{\tau \wedge t} b_i]$. 
Using monotone and dominant convergence theorem on the decomposition
\begin{equation}
    K_{t \wedge \tau} = K_\tau \id_{\{ t \ge \tau \}} + K_t \id_{\{ t < \tau \}},
\end{equation}
we have $\Exp[K_\tau \mathbbm{1}_{\tau < \infty}] \ge K_0 - \Exp[\sum_{t=1}^\tau b_t]$.
On the other hand, $K_\tau = K_{\tau-1} Y_\tau$. Conditioning on $\tau=t$, we have $K_{t-1} < 1$, so $Y_t \ge 1/K_{t-1} > 1$. Thus, using the assumption with $\beta = 1/K_{t-1}$,
\begin{align*}
\Exp[K_\tau \mathbbm{1}_{\tau < \infty}] &= \sum_{t=1}^\infty \Exp\left[ \Exp[K_{t-1} Y_t \mid \cF_{t-1}, Y_t \ge 1/K_{t-1}] \mathbbm{1}_{\tau=t} \right] \\
&\le \sum_{t=1}^\infty \Exp\left[ K_{t-1} \frac{1}{K_{t-1}} (1+\varepsilon) \mathbbm{1}_{\tau=t} \right] = (1+\varepsilon) \Pr(\tau < \infty).
\end{align*}
Combining the two inequalities yields the desired bound.
\end{proof}

\underline{Step 2: Constructing a martingale game instance.} Fix $\varepsilon,\delta > 0$. We construct a specific instance of the borrowed betting game where Casino offers two-outcome fair bets on biased coins; that is, $B_t | \cF_{t-1}$ is a distribution supported on some $\{ x,y \} \subseteq [0,\infty)$ with mean 1. 
Statistician uses the following strategy. 
\begin{itemize}
\item When $W_{t-1} > 0$, Statistician borrows $\beta_t = 0$ and bets on a biased coin with outcome $B_t \in \{ 0, 1+\varepsilon \}$. That is, $W_t \in \{0, W_{t-1} (1+ \varepsilon)\}$.
\item When $W_{t-1} = 0$, Statistician borrows $\beta_t = \delta$ and also bets on a biased coin with outcome $B_t \in \{ 0, 1+\varepsilon \}$. That is, $W_t \in \{0, \delta (1+ \varepsilon)\}$.
\end{itemize}
We call these two scenarios ``non-borrowing $t$" and ``borrowing $t$'' respectively.
Let us verify the conditions of \cref{lem:almost_sub_anti} for the process $\{K_t\}$, which we previously defined as 
\begin{equation}
  K_t = 1- \frac{h(L_t) - N_t }{h(L_t) + L_t} s(L_t) =  1- \frac{g(L_t) - W_t }{g(L_t)} s(L_t)  ,\quad
  s(L) = e^{-\int_L^\infty \frac{1}{g(\ell)} d \ell}.
\end{equation}
First, at borrowing $t$, the condition $\Exp[K_t \mid \cF_{t-1}] \ge K_{t-1} - b_t$ holds with
\begin{multline}
    K_{t-1} - \Exph[K_t \mid \cF_{t-1}] = \left(1- \frac{g(L_{t-1}) - 0 }{g(L_{t-1})} s(L_{t-1}) \right)-  \left(1- \frac{g(L_{t-1}+\delta) - \delta }{g(L_{t-1}+\delta)} s(L_{t-1} + \delta)\right)
         \\
    = s(L_{t-1}+\delta) - s(L_{t-1}) - \frac{s(L_{t-1}+\delta) }{g(L_{t-1}+\delta)}\delta =: b_t \ge 0;
\end{multline}
whereas at non-borrowing $t$, $\Exp[K_t \mid \cF_{t-1}] = K_{t-1} \ge K_{t-1} - b_t$ also holds.

Second, since borrowing $t$ must happen infinitely often, $L_t \to \infty$, $K_t \to 0$. 

Finally, let us show that
the multiplicative jumps $Y_t = K_t/K_{t-1}$ are bounded. When $W_{t-1} > 0$,
\begin{equation}
    Y_t = \frac{1- \frac{g(L_{t-1}) - W_t }{g(L_{t-1})} s(L_{t-1})}{1- \frac{g(L_{t-1}) - W_{t-1} }{g(L_{t-1})} s(L_{t-1})} \le \frac{1- \frac{g(L_{t-1}) - W_{t-1}(1+\varepsilon)}{g(L_{t-1})} s(L_{t-1})}{1- \frac{g(L_{t-1}) - W_{t-1} }{g(L_{t-1})} s(L_{t-1})} \le  \frac{1- \frac{g(L_{t-1}) - W_{t-1}(1+\varepsilon)}{g(L_{t-1})} \cdot 1}{1- \frac{g(L_{t-1}) - W_{t-1} }{g(L_{t-1})} \cdot 1}  = 1 + \varepsilon.
\end{equation}
When $W_{t-1} = 0$, 
\begin{multline}
     Y_t = \frac{1- \frac{g(L_{t}) - W_t }{g(L_{t})} s(L_{t})}{1-  s(L_{t-1})} \le  \frac{1- \frac{g(L_{t}) - \delta(1+\varepsilon) }{g(L_{t})} s(L_{t})}{1-  s(L_{t-1})}  \\
  \stackrel{[!]}{\le}  \frac{1- \frac{g(L_{t})  }{g(L_{t})} s(L_{t}) + \{ s(L_{t-1}+\delta) - s(L_{t-1}) + \varepsilon(1 - s(L_{t-1}))\}}{1-  s(L_{t-1})} = 1+\varepsilon,
\end{multline}
where the inequality marked with $[!]$ is based on the following upper bounds on $\delta \frac{s(L+\delta)}{g(L+\delta)}$: we note that
\begin{equation}
    \log \frac{s(y)}{s(x)} = \int_x^y  \frac{\d u}{g(u)} \le \log \frac{g(x) + y - x}{g(x)} \le  \log \frac{g(y)}{g(x)},
\end{equation}
meaning
$s(\ell)/g(\ell)$ is decreasing, implying $\delta \frac{s(L+\delta)}{g(L+\delta)} \le \int_L^{L+\delta} \frac{s(\ell)}{g(\ell)}\d\ell = s(L+\delta) - s(L) \le  1 - s(L)$.

Thus, $Y_t \le 1+\varepsilon$ almost surely, and the condition ``$\Exph[Y_t \mid \cF_{t-1}, Y_t \ge \beta] \le \beta(1+\varepsilon)$ for all $\beta \ge 1$'' is satisfied.

\underline{Step 3: Evaluating the lower bound.} Recall again that $s'(\ell) = s(\ell)/g(\ell)$. The total borrowing in the instance above is controlled by a Riemann integral approximation error. Let $T$ be the stopping time when $K$ borrowing $t$'s have happened.
\begin{equation}
    \sum_{t=1}^T b_t =  \sum_{k=1}^K s(k\delta) - s((k-1)\delta) - {\delta}s'(k\delta) \le \int_0^\infty s'(\ell) \d \ell - \delta\sum_{k=1}^\infty s'(k\delta).
\end{equation}
This approximation error is upper bounded by some $\mathcal{E}(\delta)$ where $\mathcal{E}(\delta)\to 0$ as $\delta \to 0$. Therefore, invoking \cref{lem:almost_sub_anti},
\[
\Pr\left( \sup_{t\ge 0} K_t \ge 1 \right) \ge \frac{K_0 - \mathcal{E}(\delta) }{1+\varepsilon}.
\]
Notice that $\Ps(\sup_t K_t \ge 1)$ corresponds exactly to $\Ps(\exists t: M_t \ge g(L_t))$.
By taking the limit $\delta, \varepsilon \to 0$,
\[
\sup \Ps(\exists t: M_t \ge g(L_t)) \ge K_0 = 1 - \frac{g(0)-1}{g(0)}\exp\left(-\int_0^\infty \frac{1}{g(\ell)}d\ell\right).
\]
This demonstrates the tightness of the bound.
\end{proof}

\end{document}

%% file: macros.tex
\renewcommand{\le}{\leqslant}

\renewcommand{\geq}{\geqslant}
\renewcommand{\ge}{\geqslant}

\newcommand{\id}{\mathbbmss 1}

\newcommand {\matl}{\left[ \begin{matrix}}
\newcommand {\matr}{\end{matrix}\right]}
\newcommand {\Exp}{ \mathbb E }

\renewcommand {\Pr}{ \mathbb P }

\renewcommand{\d}{\mathrm{d}}

\newcommand{\cA}{\mathcal{A}}

\newcommand{\cF}{\mathcal{F}}

\DeclareMathAlphabet{\mathbbmsl}{U}{bbm}{m}{sl}

\newcommand{\hwcmt}[1]{}

\newcommand{\Exph}{\Exp_{\mathsf M}}
\newcommand{\Ph}{\Pr_{\mathsf M}}

\newcommand{\Exps}{\Exp_{\mathsf S}}
\newcommand{\Ps}{\Pr_{\mathsf S}}